\documentclass[10pt]{article}

\usepackage{amsmath,amssymb,diagrams,newcent}
\usepackage{a4}

\newtheorem{theorem}{Theorem}

\newcommand{\np}{{\par \noindent}}

\newcommand{\C}{\mathbb{C}}

\newcommand{\N}{\mathbb{N}}

\newcommand{\Z}{\mathbb{Z}}
\newcommand{\PP}{\mathbb{P}}

\newcommand{\atn}{\text{\bf @}_{\text{\bf n}}}

\newcommand{\wisk}[1]{{\text{\bf \usefont{OT1}{pag}{m}{n} #1}}}

\begin{document}

\title{Noncommutative compact manifolds constructed from quivers}
\author{Lieven Le Bruyn\thanks{Research director of the FWO (Belgium)} 
\\ Universitaire Instelling Antwerpen \\
B-2610 Antwerp (Belgium)  \\
{\tt lebruyn$@\mbox{wins.uia.ac.be}$}}

\date{\empty}
\maketitle

\begin{abstract}
The moduli spaces of $\theta$-semistable representations of a finite quiver can be packaged
together to form a noncommutative compact manifold.
\end{abstract}

\np
If noncommutative affine schemes are geometric objects associated to affine associative $\C$-algebras,
affine smooth noncommutative varieties ought to correspond to
{\it quasi-free (or formally smooth) algebras} (having the lifting property for
algebra morphisms modulo nilpotent ideals). Indeed, J. Cuntz and D. Quillen have shown that for an algebra to
have a rich theory of differential forms allowing natural connections it must be quasi-free
\cite[Prop. 8.5]{CuntzQuillen}.

\np
M. Kontsevich and A. Rosenberg introduced {\it noncommutative spaces}
generalizing the notion of stacks to the noncommutative case \cite[\S 2]{KontRos}. It is hard to
construct noncommutative compact manifolds in this framework, due to the scarcity of
faithfully flat extensions for quasi-free algebras. An alternative was
outlined by M. Kontsevich in \cite{Kontold} and made explicit in \cite[\S 1]{KontRos}
(see also \cite{LBngatn} and \cite{Kontnew}). Here, the geometric object corresponding to the
quasi-free algebra $A$ is the collection $(\wisk{rep}_n A)_n$ where $\wisk{rep}_n A$ is the affine
$GL_n$-scheme of $n$-dimensional representations of $A$. As $A$ is quasi-free each $\wisk{rep}_n A$
is smooth and endowed with Kapronov's formal noncommutative structure \cite{Kapranov}.
Moreover, this collection has equivariant {\it sum-maps}
$\wisk{rep}_n A \times \wisk{rep}_m A \rTo \wisk{rep}_{m+n} A$.

\np
We define a {\it noncommutative compact manifold} to be a collection $(Y_n)_n$ of projective varieties
such that $Y_n$ is the quotient-scheme of a smooth $GL_n$-scheme $X_n$ which is locally isomorphic to 
$\wisk{rep}_n A_{\alpha}$ for a fixed set of quasi-free algebras $A_{\alpha}$, is endowed
 with a formal noncommutative
structure and there are equivariant sum-maps $X_m \times X_n \rTo X_{m+n}$. In this note we will 
construct of a large class of examples.

\np
An illustrative example : let $M_{\PP_2}(n;0,n)$ be the moduli space of semi-stable vectorbundles
of rank $n$ over the projective plane $\PP_2$ with Chern-numbers $c_1=0$ and $c_2=n$, then the
collection $(M_{\PP_2}(n;0,n))_n$ is a noncommutative compact manifold. In general, let $Q$ be a
quiver on $k$ vertices {\it without oriented cycles} and let $\theta = (\theta_1,\hdots,\theta_k) 
\in \Z^k$. For a finite dimensional representation $N$ of $Q$ with dimension vector $\alpha =
(a_1,\hdots,a_k)$ we denote $\theta(N) = \sum_i \theta_i a_i$ and $d(\alpha) = \sum_i a_i$. A representation $M$ of $Q$ is called
$\theta$-{\it semistable} if $\theta(M)=0$ and $\theta(N) \geq 0$ for every subrepresentation $N$ of $M$.
A. King studied the moduli spaces $M_Q(\alpha, \theta)$ of $\theta$-semistable representations of $Q$
of dimension vector $\alpha$ and proved that these are projective varieties \cite[Prop 4.3]{King}. 
We will prove the following result.

\begin{theorem} With notations as above, the collection of projective varieties
\[
( \bigsqcup_{d(\alpha) = n} M_Q(\alpha,\theta)~)_n \]
is a noncommutative compact manifold.
\end{theorem}

\np
The claim about moduli spaces of vectorbundles on $\PP_2$ follows by considering the quiver
$\bullet \pile{\rTo \\ \rTo \\ \rTo} \bullet$ and $\theta = (-1,1)$.

\np
Let $C$ be a smooth projective curve of genus $g$ and $M_C(n,0)$ the moduli space of semi-stable
vectorbundles of rank $n$ and degree $0$ over $C$. We expect the collection $(M_C(n,0))_n$ to be a
noncommutative compact manifold.

\section{The setting.}

\np
Let $Q$ be a {\it quiver} on a finite set $Q_v = \{ v_1,\hdots,v_k \}$ of
vertices having a finite set $Q_a$ of arrows.
We assume that $Q$ has {\it no oriented cycles}. 

\np
The path algebra $\C Q$ has as underlying $\C$-vectorspace basis the set of all oriented paths
in $Q$, including those of length zero which give idempotents corresponding to the vertices $v_i$. Multiplication in
$\C Q$ is induced by (left) concatenation of paths. $\C Q$ is a finite dimensional quasi-free
algebra. 

\np
Let $\alpha = (a_1,\hdots,a_k)$ be a {\it dimension vector} such that $d(\alpha) = n$.
Let $\wisk{rep}_{Q}(\alpha)$ be the affine space of $\alpha$-dimensional representations of the quiver
$Q$. That is,
\[
\wisk{rep}_{Q}(\alpha) = \bigoplus_{\underset{j}{\bullet} \lTo^a \underset{i}{\bullet}}
M_{a_j \times a_i}(\C) \]
$GL(\alpha) = GL_{a_1} \times \hdots \times GL_{a_k}$ acts on this space via basechange in the 
vertexspaces. For $\theta = (\theta_1,\hdots,\theta_k) \in \Z^k$ we denote with
$\wisk{rep}_Q^{ss}(\alpha,\theta)$ the open (possibly empty) subvariety of $\theta$-semistable
representations in $\wisk{rep}_Q(\alpha)$. Applying results of A. Schofield \cite{schofpaper} there
is an algorithm to determine the $(\alpha,\theta)$ such that $\wisk{rep}_Q^{ss}(\alpha,\theta) \not= 
\emptyset$. Consider the diagonal embedding of $GL(\alpha)$ in $GL_n$ and the quotient morphism
\[
X_n = \bigsqcup_{d(\alpha) = n} GL_n \times^{GL(\alpha)} rep^{ss}_Q(\alpha,\theta) \rOnto^{\pi_n} Y_n = \bigsqcup_{d(\alpha) = n}
M_Q(\alpha,\theta) . \]
Clearly, $X_n$ is a smooth $GL_n$-scheme and the direct sum of representations induces sum-maps
$X_m \times X_n \rTo X_{m+n}$ which are equivariant with respect to $GL_m \times GL_n \rInto GL_{m+n}$.
$Y_n$ is a projective variety by \cite[Prop. 4.3]{King} and its points correspond to isoclasses of
$n$-dimensional representations of $\C Q$ which are direct sums of $\theta$-stable representations by
\cite[Prop. 3.2]{King}. Recall that a $\theta$-semistable representation $M$ is called $\theta$-{\it stable}
provided the only subrepresentations $N$ with $\theta(N)=0$ are $M$ and $0$.

\section{Universal localizations.}

\np
We recall the notion of {\it universal
localization} and refer to \cite[Chp. 4]{Schofield} for full details.
Let $A$ be a $\C$-algebra and $\wisk{projmod}~A$ the category of finitely generated projective left $A$-modules.
Let $\Sigma$ be some class of maps in
this category. In \cite[Chp. 4]{Schofield}
it is shown that there exists an algebra map $A \rTo^{j_{\Sigma}} A_{\Sigma}$ with the universal
property that the maps $A_{\Sigma} \otimes_A \sigma$ have an inverse for all $\sigma \in \Sigma$.
$A_{\Sigma}$ is called the universal localization of $A$ with respect to the set of maps $\Sigma$.
In the special case when $A$ is the path algebra $\C Q$ of a quiver on $k$ vertices,  we can
identify the isomorphism classes in $\wisk{projmod}~\C Q$ with $\N^k$. To each vertex $v_i$ corresponds
an {\it indecomposable} projective left $\C Q$-ideal $P_i$ having as $\C$-vectorspace basis all paths
in $Q$ {\it starting at} $v_i$. We can also determine the space of homomorphisms
\[
Hom_{\C Q}(P_i,P_j) = \bigoplus_{\underset{i}{\bullet} \lDotsto^p \underset{j}{\bullet}} \C p \]
where $p$ is an oriented path in $Q$ starting at $v_j$ and ending at $v_i$. Therefore, any
$A$-module morphism $\sigma$ between two projective left modules
\[
P_{i_1} \oplus \hdots \oplus P_{i_u} \rTo^{\sigma} P_{j_1} \oplus \hdots \oplus P_{j_v} \]
can be represented by an $u \times v$ matrix $M_{\sigma}$ whose $(p,q)$-entry $m_{pq}$ is a linear
combination of oriented paths in $Q$ starting at $v_{j_q}$ and ending at $v_{i_p}$.

\np
Now, form an $v \times u$ matrix $N_{\sigma}$ of free variables $y_{pq}$ and consider the algebra
$\C Q_{\sigma}$ which is the quotient of the free product $\C Q \ast \C \langle y_{11},\hdots,y_{uv}
\rangle$ modulo the ideal of relations determined by the matrix equations
\[
M_{\sigma}. N_{\sigma} = \begin{bmatrix} v_{i_1} & & 0 \\
& \ddots & \\
0 & & v_{i_u} \end{bmatrix} \quad \quad 
N_{\sigma} . M_{\sigma} = \begin{bmatrix}
v_{j_1} & & 0 \\ & \ddots & \\
0 & & v_{j_v} \end{bmatrix}. \]

\np
Repeating this procedure for every $\sigma \in \Sigma$ we obtain the universal localization
$\C Q_{\Sigma}$. Observe that if $\Sigma$ is a finite set of maps, then the universal localization
$\C Q_{\Sigma}$ is an affine algebra.

\np
It is easy to see that $\C Q_{\Sigma}$ is quasi-free and that
the representation space $\wisk{rep}_n~\C Q_{\sigma}$ is an
open subscheme (but possibly empty) of $\wisk{rep}_n~\C Q$. Indeed, if $m = (m_a)_a \in \wisk{rep}_{Q}(\alpha)$,
then $m$ determines a point in $\wisk{rep}_n~\C Q_{\Sigma}$ if and only if the matrices
$M_{\sigma}(m)$ in which the arrows are all replaced by the matrices $m_a$ are invertible for all
$\sigma \in \Sigma$. In particular, this induces numerical conditions on the dimension
vectors $\alpha$ such that $\wisk{rep}_n \C Q_{\Sigma} \not= \emptyset$. Let $\alpha = (a_1,\hdots,a_k)$
be a dimension vector such that $\sum a_i = n$ then every $\sigma \in \Sigma$ say with
\[
P_1^{\oplus e_1} \oplus \hdots \oplus P_k^{\oplus e_k} \rTo^{\sigma} P_1^{\oplus f_1} \oplus
\hdots \oplus P_k^{\oplus f_k} \]
gives the numerical condition
$
e_1 a_1 + \hdots + e_k a_k = f_1 a_1 + \hdots + f_k a_k$.

\section{Local structure.}

\np
Fix $\theta = (\theta_1,\hdots,\theta_k) \in \Z^k$ and let $\Sigma
= \cup_{z \in \N_+} \Sigma_z$ where $\Sigma_z$ is the set of all morphisms $\sigma$
\[
P_{i_1}^{\oplus z\theta_{i_1}} \oplus \hdots \oplus P_{i_u}^{\oplus z\theta_{i_u}} \rTo^{\sigma}
P_{j_1}^{\oplus - z\theta_{j_1}} \oplus \hdots \oplus P_{j_v}^{\oplus - z\theta_{j_v}} \]
where $\{ i_1,\hdots,i_u \}$ (resp. $\{ j_1,\hdots,j_v \}$) is the set of indices $1 \leq i \leq k$
such that $\theta_i > 0$ (resp. $\theta_i < 0$). Fix a dimension vector
$\alpha$ with $\langle \theta,\alpha \rangle = 0$, then $\theta$ determines a character $\chi_{\theta}$ on
$GL(\alpha)$ defined by $\chi_{\theta}(g_1,\hdots,g_k) = \prod det(g_i)^{\theta_i}$. With notations as
before, the function $d_{\sigma}(m) = det(M_{\sigma}(m))$ for $m \in \wisk{rep}_Q(\alpha)$ is a {\it semi-invariant}
of weight $z\chi_{\theta}$ in $\C[\wisk{rep}_Q(\alpha)]$ if $\sigma \in \Sigma_z$.

\np
The open subset $X_{\sigma}(\alpha) = \{ m \in \wisk{rep}_Q(\alpha) \mid d_{\sigma}(m) \not= 0 \}$ 
consists of $\theta$-semistable representations which are also $n$-dimensional representations of the
universal localization $\C Q_{\sigma}$. Under this correspondence $\theta$-stable representations
correspond to simple representations of $\C Q_{\sigma}$. If we denote
\[
X_{\sigma,n} = \bigsqcup_{d(\alpha) = n} GL_n \times^{GL(\alpha)} X_{\sigma}(\alpha) \rInto X_n \]
then $X_{\sigma, n} = \wisk{rep}_n \C Q_{\sigma}$ and the restriction of $\pi_n$ to $X_{\sigma,n}$
is the $GL_n$-quotient map $\wisk{rep}_n \C Q_{\sigma} \rOnto \wisk{fac}_n \C Q_{\sigma}$ which sends
an $n$-dimensional representation to the isomorphism class of the semi-simple $n$-dimensional representation
of $\C Q_{\sigma}$ given by the sum of the Jordan-H\"older components, see \cite[2.3]{LBngatn}. As the
semi-invariants $d_{\sigma}$ for $\sigma \in \Sigma$ cover the moduli spaces $M_Q(\alpha,\theta)$ this
proves the local isomorphism condition for the collection $(Y_n)_n$.

\np
A point $y \in Y_n$ determines a unique closed orbit in $X_n$ corresponding to a representation
\[
M_y = M_1^{\oplus e_1} \oplus \hdots \oplus M_l^{\oplus e_l} \]
with the $M_i$ $\theta$-stable representations occurring in $M_y$ with multiplicity $e_i$. The local
structure of $Y_n$ near $y$ is completely determined by a {\it local quiver} $\Gamma_y$ on $l$
vertices which
usually
has loops and oriented cycles and a dimension vector $\beta_y = (e_1,\hdots,e_l)$. The quiver-data
 $(\Gamma_y,\beta_y)$ is determined by the canonical $A_{\infty}$-structure on 
 $Ext^{\ast}_{\C Q}(M_y,M_y)$. As $\C Q$ is quasi-free, this ext-algebra has only components in degree zero
 (determining the vertices and the dimension vector $\beta_y$) and degree one (giving the
 arrows in $\Gamma_y$). 
 
\np
Using \cite[Thm 4.7]{Schofield} and the correspondence between $\theta$-stable representations and
simples of universal localizations, the local structure is the one outlined in \cite[2.5]{LBngatn}.
In particular, it can be used to locate the singularities of the projective varieties $Y_n$.

\section{Formal structure.}

\np
In \cite{Kapranov} M. Kapranov computes the formal neighborhood of commutative manifolds embedded
in noncommutative manifolds.
Equip a $\C$-algebra $R$ with the {\it commutator filtration} having as part of degree $-d$
\[
F_{-d} = \sum_{m} \sum_{i_1 + \hdots i_m = d} RR^{Lie}_{i_1}R \hdots R R^{Lie}_{i_m}R \]
where $R^{Lie}_i$ is the subspace spanned by all expressions 
$[r_1,[r_2,[ \hdots,[r_{i-1},r_i] \hdots]$ containing $i-1$ instances of Lie brackets. 
We require that for $R_{ab} = \tfrac{R}{F_{-1}}$ affine smooth, the algebras $\tfrac{R}{F_{-d}}$ have the lifting
property modulo nilpotent algebras in the category of $d$-nilpotent algebras (that is, those
such that $F_{-d} = 0$). The micro-local structuresheaf with respect to the commutator
filtration then defines a sheaf of noncommutative algebras on $\wisk{spec} R_{ab}$, the
{\it formal structure}. Kapranov shows that in the affine case there exists an essentially
unique such structure. For arbitrary manifolds there is an obstruction to the existence of
a formal structure and when it exists it is no longer unique. We refer to \cite[4.6]{Kapranov} for
an operadic interpretation of these obstructions.

\np
We will write down the formal structure on the affine open subscheme $\wisk{rep}_n \C Q_{\Gamma}$
of $X_n$ where $\Gamma$ is a finite subset of $\Sigma$. Functoriality of this construction then
implies that one can glue these structures together to define a formal structure on $X_n$ finishing
the proof of theorem 1.

\np
If $A$ is an affine quasi-free algebra, the formal structure on $\wisk{rep}_n A$ is given by the micro-structuresheaf
for the commutator filtration on the affine algebra
\[
\sqrt[n]{A} = A \ast M_n(\C)^{M_n(\C)} = \{ p \in A \ast M_n(\C) \mid p.(1 \ast m) = 
(1 \ast m).p \ \forall m \in M_n(\C) \}
\]
This follows from the fact that $\sqrt[n]{A}$ is again quasi-free by the coproduct theorems,
\cite[\S 2]{Schofield}. Specialize to the case when $A = \sqrt[n]{\C Q_{\Gamma}}$. Consider the
extended quiver $\hat{Q}(n)$ by adding one vertex $v_0$ and for every vertex $v_i$ in $Q$ we add
$n$ arrows from $v_0$ to $v_i$ denoted $\{ x_{i1},\hdots,x_{in} \}$. 
Consider the morphism between projective left $\C \hat{Q}(n)$-modules
\[
P_1 \oplus P_2 \oplus \hdots \oplus P_k \rTo^{\tau} \underbrace{P_0 \oplus \hdots \oplus P_0}_n \]
determined by the matrix
\[
M_{\tau} = \begin{bmatrix} x_{11} & \hdots & \hdots & x_{1n} \\
\vdots & & & \vdots \\
x_{k1} & \hdots & \hdots & x_{kn} \end{bmatrix}. \]
Consider the universal localization $B = \C \hat{Q}(n)_{\Gamma \cup \{ \tau \} }$. Then,
$\sqrt[n]{\C Q_{\Gamma}} = v_0 B v_0$ the algebra of oriented loops based at $v_0$.

\section{Odds and ends.}

\np
One can build a global combinatorial object from the universal localizations $\C Q_{\Gamma}$ with
$\Gamma$ a finite subset of $\Sigma$ and gluings coming from unions of these sets. This example may
be useful to modify the Kontsevich-Rosenberg proposal of noncommutative spaces to the quasi-free
world.

\np
Finally, allowing oriented cycles in the quiver $Q$ one can repeat the foregoing verbatim and obtain a
projective space bundle over the collection $(\wisk{fac}_n \C Q)_n$.

\end{document}